\documentclass[12pt]{article}
\usepackage{amssymb}
\usepackage{amsmath}
\usepackage{indentfirst,latexsym}
\usepackage{mathrsfs}
\usepackage{graphicx}
\usepackage{fancyhdr}

\date{}
\begin{document}

\title{Occurrences of consecutive patterns of length $3$ in $3\textrm{-}1\textrm{-}2$-avoiding permutations}
\author{Marilena Barnabei, Flavio Bonetti, and Matteo Silimbani \thanks{
Dipartimento di Matematica - Universit\`a di Bologna}} \maketitle

\noindent {\bf Abstract} We exploit Krattenthaler's bijection
between the set $S_n(3\textrm{-}1\textrm{-}2)$ of permutations in
$S_n$ avoiding the classical pattern $3\textrm{-}1\textrm{-}2$ and
Dyck $n$-paths to study the distribution of every consecutive
pattern of length $3$ on the set $S_n(3\textrm{-}1\textrm{-}2)$.
We show that these consecutive patterns split into $3$
equidistribution classes, by means of an involution on Dyck paths
due to E.Deutsch. In addition, we state equidistribution theorems concerning triplets of
statistics relative to the occurrences of the consecutive patterns of length $3$ in a permutation.\\


\section{Introduction}

\noindent Let $\sigma\in S_n$ and $\tau\in S_k$, $k\leq n$, be two
permutations. We say that $\sigma$ \emph{contains the pattern}
$\tau$ if $\sigma$ contains a subsequence order-isomorphic to
$\tau$. We say that $\sigma$ \emph{avoids} $\tau$ if such a
subsequence does not exist. The subject of pattern avoiding
permutations was initiated by Simion and Schimdt \cite{ss}, and,
after that, a large literature on this topic has blossomed.
Problems treated include counting permutations avoiding a pattern
or a set of patterns, or containing patterns a specified number of
times.

\noindent More recently, Babson and Steingr\'imsson \cite{bab}
introduced generalized permutation patterns where two adjacent
letters in a pattern may be required to be adjacent in the
permutation. A number of interesting results on generalized
patterns were obtained by several authors in recent years (for an
extensive survey, see \cite{steim}).

\noindent One particular case of generalized patterns are
consecutive patterns. For a subsequence of a permutation to be an
occurrence of a consecutive pattern, its elements have to appear
in adjacent positions of the permutation. We write a ''classical''
pattern with dashes between symbols, while a consecutive pattern
will be written without dashes, accordingly with the most common
notation (see \cite{steim}). A number of results for the
enumeration of permutations by consecutive patterns have recently
been obtained (see, e.g., \cite{avg}, \cite{eli}, \cite{sergi},
\cite{kms}, \cite{mend}, \cite{raw}, \cite{warl}, and
\cite{warlii}).

\noindent Many well known integer sequences arise in enumerative
problems concerning permutations avoiding a pattern $\tau$ or
containing it a fixed number of times. In particular, in \cite{ss}
it has been shown that the number of permutations avoiding any
classical pattern $\tau\in S_3$ equals the $n$-th Catalan number
$C_n$. Afterward, Krattenthaler \cite{kratt} described a bijection
between Dyck paths of semilength $n$ and the set
$S_n(3\textrm{-}1\textrm{-}2)$ of permutations avoiding
$3\textrm{-}1\textrm{-}2$.\\

\noindent In this paper, we study the distribution of the five
non-trivial consecutive patterns of length $3$ on the set of
$3\textrm{-}1\textrm{-}2$-avoiding permutations.

\noindent More precisely, for every consecutive pattern $\tau$ of
length $3$, $\tau\neq 312$, we study the bivariate generating
function
$$A^{\tau}(t,z)=\sum_{n,k\geq 0} a^{\tau}_{n,k}t^n z^k,$$
where $a_{n,k}^{\tau}$ is the number of permutations in
$S_n(3\textrm{-}1\textrm{-}2)$ containing $k$ occurrences of the
consecutive pattern $\tau$. First of all, we prove that each
occurrence of such a $\tau$ in a permutation $\sigma$ correspond
bijectively to a peculiar configuration in the Dyck path
associated to $\sigma$ by Krattenthaler's bijection. This
correspondence allows us to show that the five non-trivial
patterns of length $3$ split into $3$ classes (i.e., $\{213\}$,
$\{123,321\}$, and $\{132,231\}$), so that two patterns in the
same class are equidistributed on $S_n(3\textrm{-}1\textrm{-}2)$,
namely, the corresponding bivariate generating functions coincide.
These equidistribution results are obtained applying an involution
on Dyck paths due to Deutsch \cite{deu}.

\noindent For each one of these classes we can choose a
representative such that the distribution of the corresponding
Dyck path configuration has been determined (see \cite{deu2} and
\cite{sapo}). Hence, we get the bivariate generating function of
the distribution of each one of the five non-trivial consecutive
patterns. This allows us to get an explicit expression for the
coefficients of the power series $A^{\tau}(t,z)$ for every
consecutive pattern $\tau$ of length $3$.

\noindent As a fallout, we obtain a formula for the number of
permutations in $S_n(3\textrm{-}1\textrm{-}2)$ that avoid any
consecutive pattern $\tau\in S_3$. In the two cases $\tau=321$ and
$\tau=123$, we get the Motzkin numbers, and present two bijections
$\nu:S_n(3\textrm{-}1\textrm{-}2,321)\to\mathscr{M}_n$ and
$\mu:S_n(3\textrm{-}1\textrm{-}2,123)\to\mathscr{M}_n$, where
$\mathscr{M}_n$ denotes the set of Motzkin $n$-paths.

\noindent In the last section, we show how previous results yield
equidistribution theorems concerning triplets of statistics that
associate the number of occurrences of a given consecutive pattern
with each permutation in $S_n(3\textrm{-}1\textrm{-}2)$.

\section{Preliminaries}
\label{bestia}

\subsection{Lattice paths}

\noindent A \emph{Dyck path} of semilength $n$ (or Dyck $n$-path)
is a lattice path starting at $(0,0)$, ending at $(2n,0)$, and
never going below the $x$-axis, consisting of up steps $U=(1,1)$
and down steps $D=(1,-1)$. A \emph{Motzkin path} of length $n$ (or
Motzkin $n$-path) is a lattice path starting at $(0,0)$, ending at
$(n,0)$, and never going below the $x$-axis, consisting of up
steps $U=(1,1)$, horizontal steps $H=(1,0)$, and down steps
$D=(1,-1)$.

\noindent Dyck paths of semilength $n$ are counted by the $n$-th
Catalan number $C_n$, while Motzkin paths of length $n$ are
counted
by the $n$-th Motzkin number $M_n$.\\

\noindent A Dyck path can be regarded as a word over the alphabet
$\{U,D\}$ such that any prefix contains at least as many symbols
$U$ as symbols $D$. A \emph{subword} of a Dyck path $P$ is a
subsequence of consecutive steps in $P$.

\noindent An \emph{irreducible} Dyck path is a Dyck path that does
not touch the $x$-axis except for the origin and the final
destination. An \emph{irreducible component} of a Dyck path $P$ is
a maximal irreducible Dyck subpath of $P$.

\noindent We list some notions on Dyck paths that will be used in
the following. A \emph{run} (respectively \emph{fall}) of a Dyck
path is a maximal subword consisting of up (resp. down) steps. A
\emph{return} of a Dyck path is a down step landing on
the $x$-axis.\\

\begin{figure}[h]
\begin{center}
\includegraphics[bb=70 453 331 708,width=.6\textwidth]{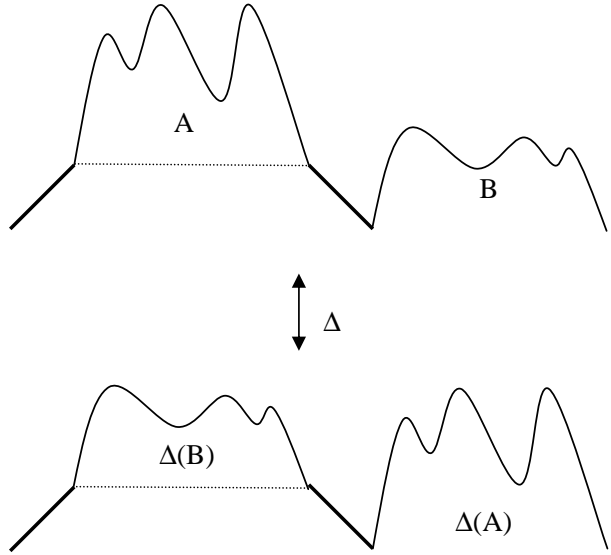}
\caption{The involution $\Delta$.}\label{nonusata}
\end{center}
\end{figure}

\begin{figure}[h]
\begin{center}
\includegraphics[bb=74 392 513 703,width=.5\textwidth]{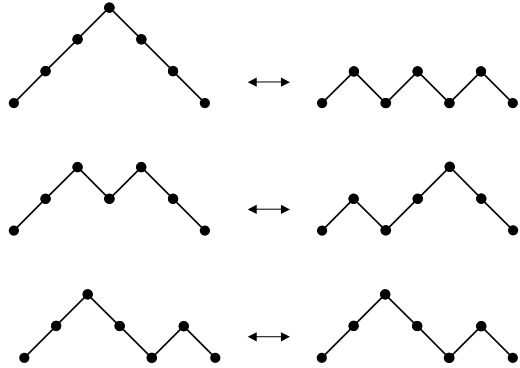}
\caption{The involution $\Delta$ on Dyck paths of semilength
$3$.}\label{piccoli}
\end{center}
\end{figure}

\noindent We now describe an involution $\Delta$ on Dyck paths due
to Deutsch \cite{deu}. Consider a Dyck path $P$ and decompose it
according to its first return as $P=U\,A\,D\,B$, where $A$ and $B$
are (possibly empty) Dyck paths. Then, the Dyck path $\Delta(P)$
is recursively determined by the following rules (see Figure
\ref{nonusata}):
\begin{itemize}
\item if $P$ is empty, so is $\Delta(P)$;
\item otherwise, $\Delta(P)=U\,\Delta(B)\,D\,\Delta(A)$.
\end{itemize}
It is easily checked that the two paths $P$ and $\Delta(P)$ have
the same length. In Figure \ref{piccoli} we show how the
involution $\Delta$ acts on Dyck paths of semilength $3$.

\subsection{Restricted permutations}

\noindent Let $\sigma\in S_n$ and $\pi\in S_k$, $k\leq n$, be two
permutations. The permutation $\sigma$ \emph{contains the pattern}
$\pi$ if there exists a subsequence
$\sigma(i_1),\sigma(i_2),\ldots,\sigma(i_k)$ with $1\leq
i_1<i_2<\cdots<i_k\leq n$ that is order-isomorphic to $\pi$.
Moreover, if $i_1,i_2,\ldots,i_k$ are consecutive integers - i.e.
$\sigma(i_1),\sigma(i_2),\ldots,\sigma(i_k)$ is a subword of
$\sigma$ -, then we say that $\sigma$ contains the
\emph{consecutive pattern} $\pi$. In order to avoid confusion, we
write a ''classical'' pattern with dashes between symbols, while
 a consecutive pattern will be written without dashes.

\noindent For example, the permutation
$\sigma=4\,3\,1\,7\,2\,5\,6$ contains $5$ occurrences of the
pattern $3\textrm{-}1\textrm{-}2$ (namely,
$412$,$312$,$725$,$726$,$756$), but only one
occurrence of the consecutive pattern $312$ (namely, $725$).\\

\noindent The permutation $\sigma$ \emph{avoids} the (consecutive)
pattern $\pi$ if $\sigma$ does not contain $\pi$. For example, the
permutation $\sigma=5\,2\,1\,3\,4$ avoids the consecutive pattern
$312$, but it does not avoid the pattern
$3\textrm{-}1\textrm{-}2$, since it contains, for instance, the
subsequence $513$.

\noindent We denote by $S_n(\tau_1,\ldots,\tau_k)$ the set of
permutations in $S_n$ that avoid simultaneously the patterns
$\tau_1,\ldots,\tau_k$.

\noindent We say that two permutation statistics $f$ and $g$ are
\emph{equidistributed} on a set $A\subseteq S_n$, if
$$\sum_{\sigma\in A}x^{f(\sigma)}=\sum_{\sigma\in
A}x^{g(\sigma)}.$$\\

\noindent In the present paper we are interested in studying the
distribution of an arbitrary consecutive pattern of length $3$ on
the set $S_n(3\textrm{-}1\textrm{-}2)$. It is well known that the
permutation $\sigma$ avoids $3\textrm{-}1\textrm{-}2$ if and only
if it can be written as follows:
$$\sigma=m_1\,w_1\,m_2\,w_2\,\ldots\,m_k\,w_k,$$
where
\begin{itemize} \item the integers $m_i$ are the left-to-right
maxima of $\sigma$ (where a \emph{left-to-right maximum} of a
permutation $\sigma$ is an integer $\sigma(i)$ such that
$\sigma(i)>\sigma(j)$ for every $j<i$);
\item for every symbol $a$ appearing in one of the words $w_i$, consider the suffix $\sigma_a$ of $\sigma$
starting with $a$. Then, $a$ is the gratest symbol in $\sigma_a$
among those that are less then $m_i$.
\end{itemize}
\noindent Some results in this direction can be found in
\cite{mandue}, where the author exhibits the generating functions
for the number of permutations on $n$ letters avoiding
$1\textrm{-}3\textrm{-}2$ (or containing $1\textrm{-}3\textrm{-}2$
exactly once) and an arbitrary generalized pattern $\tau$ on $k$
letters, or containing $\tau$ exactly once.\\

\noindent We describe the bijection $K$ between
$3\textrm{-}1\textrm{-}2$-avoiding permutations and Dyck paths due
to Krattenthaler \cite{kratt}.

\noindent Give a permutation $\sigma\in
S_n(3\textrm{-}1\textrm{-}2)$,
$\sigma=m_1\,w_1\,m_2\,w_2\,\ldots\,m_k\,w_k$ the Dyck path
$K(\sigma)$ of semilength $n$ is constructed as follows: start
with $m_1$ up steps followed by $|w_1|+1$ down steps. Then add
$m_2-m_1$ up steps followed by $|w_2|+1$ down steps, and so on.\\

\begin{figure}[h]
\begin{center}
\includegraphics[bb=57 583 512 712,width=.7\textwidth]{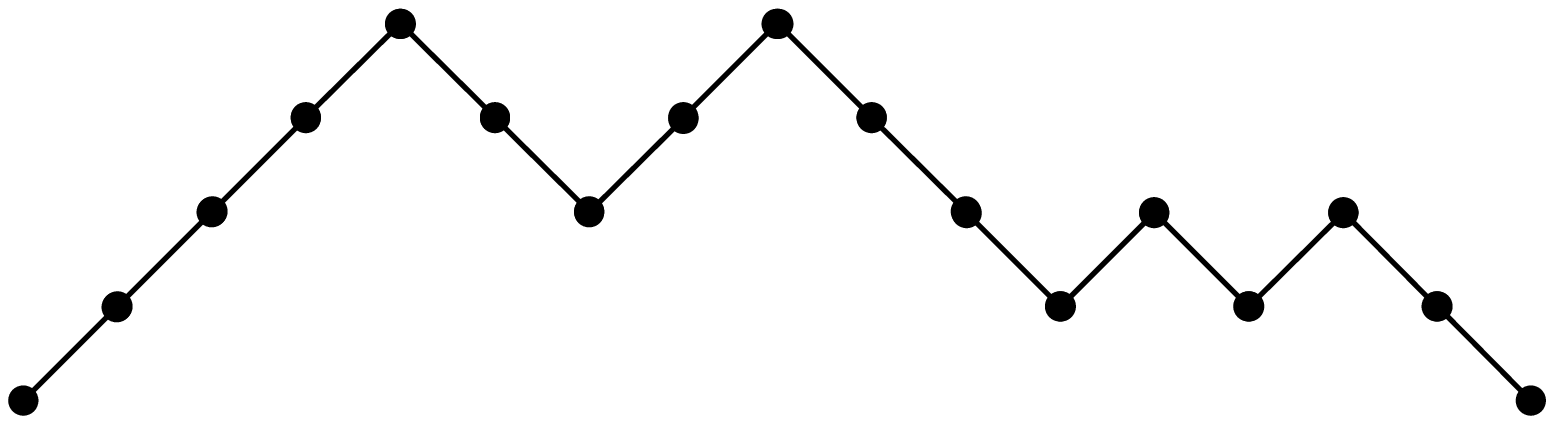}
\caption{The Dyck path $K(4\,3\,6\,5\,2\,7\,8\,1)$.}\label{kler}
\end{center}
\end{figure}

\noindent For example, the $3\textrm{-}1\textrm{-}2$-avoiding
permutation $\sigma=4\,3\,6\,5\,2\,7\,8\,1$ is mapped to the Dyck
path in Figure \ref{kler}.

\section{Consecutive patterns}

\noindent In this section  we describe the bivariate generating
function
$$A^{\tau}(t,z)=\sum_{n,k\geq 0} a^{\tau}_{n,k}t^n z^k,$$
where $a_{n,k}^{\tau}$ is the number of permutations in
$S_n(3\textrm{-}1\textrm{-}2)$ containing $k$ occurrences of the
consecutive pattern $\tau$, for every $\tau\in S_3$. Of course, we
do not consider the case $\tau=312$.

\subsection{The pattern $213$}
\label{semi}

\noindent We recall that a $3\textrm{-}1\textrm{-}2$-avoiding
permutation can be written as
$$\sigma=m_1\,w_1\,m_2\,w_2\,\ldots\,m_k\,w_k,$$ where the
integers $m_i$ are the left-to-right maxima of $\sigma$ and the
(possibly empty) subwords $w_j$ are decreasing. We note that
occurrences of the consecutive pattern $213$ in $\sigma$
correspond bijectively to nonempty subwords $w_i$, $i<k$. In fact,
if $|w_i|>0$, $i<k$, then $\sigma$ contains the $213$-subword
$b\,a\,m_{i+1}$, where $a$ is the last element in $w_i$ and $b$ is
either the second last element in $w_i$, or $b=m_i$, when
$|w_i|=1$. It is easily seen that each one of such subwords $w_i$,
in turn, corresponds to an occurrence of $DDU$ in the Dyck path
$K(\sigma)$. The distribution of the subword $DDU$ on Dyck paths
is well known (see \cite{oeis} seq.$\,$A091894). More precisely,
in \cite{deu2}, the author determines a functional equation
satisfied by the bivariate generating function of this
distribution. Hence, we deduce the following expression for the
generating function $A^{213}(t,z)$:

\newtheorem{yama}{Theorem}
\begin{yama} We have
$$A^{213}(t,z)=\frac{1-2t+2tz-\sqrt{(1-2t)^2-4t^2z}}{2tz}$$
that yields
$$a^{213}_{n,k}=2^{n-2k-1}C_k{n-1\choose 2k}$$
where $C_k$ is the $k$-th Catalan number.
\end{yama}
\begin{flushright}
$\diamond$
\end{flushright}

\noindent In particular:

\newtheorem{frt}[yama]{Proposition}
\begin{frt}
The number of permutations in $S_n$ that avoid both the pattern
$3\textrm{-}1\textrm{-}2$ and the consecutive pattern $213$ is
$$|S_n(3\textrm{-}1\textrm{-}2,213)|=2^{n-1}.$$
\end{frt}
\begin{flushright}
$\diamond$
\end{flushright}

\subsection{The pattern $321$}

\noindent Consider a $3\textrm{-}1\textrm{-}2$-avoiding
permutation
$$\sigma=m_1\,w_1\,m_2\,w_2\,\ldots\,m_k\,w_k.$$
The consecutive pattern $321$ occurs in $\sigma$ if and only if at
least one among the subwords $w_j$ has length greater than one.
More precisely, if $|w_j|=t>0$, the subword $m_j\,w_j$ contains
$t-1$ occurrences of $321$, since the elements in $w_j$ appear in
decreasing order. Note that a subword $w_j$ of length $t$
corresponds to a fall $F_j$  of length $t+1$ of the Dyck path
$K(\sigma)$. This implies that occurrences of $321$ in $m_j\,w_j$
correspond bijectively to occurrences of $DDD$ in $F_j$. The
distribution of $DDD$ on Dyck paths is well known (see \cite{oeis}
seq.$\,$A092107). In \cite{sapo}, the author deduces a functional
equation satisfied by the bivariate generating function of this
distribution. These considerations imply that:
\newtheorem{ehio}[yama]{Theorem}
\begin{ehio} We have
$$A^{321}(t,z)=\frac{1-t+tz-\sqrt{1-2t-3t^2+tz(tz+2t-2)}}{2t(t+z-tz)}$$
that yields
$$a^{321}_{n,k}=\frac{1}{n+1}\sum_{j=0}^k(-1)^{k-j}{n+j\choose n}{n+1\choose k-j}\sum_{i=j}^{\left\lfloor\frac{n+j}{2}\right\rfloor}
{n+j+1-k\choose i+1}{n-i\choose i-j}.$$
\end{ehio}
\begin{flushright}
$\diamond$
\end{flushright}
\noindent In particular:
\newtheorem{exact}[yama]{Proposition}
\begin{exact}\label{arrsub}
The number of permutations in $S_n$ that avoid both the pattern
$3\textrm{-}1\textrm{-}2$ and the consecutive pattern $321$ is
$$|S_n(3\textrm{-}1\textrm{-}2,321)|=M_n,$$
where $M_n$ is the $n$-th Motzkin number.
\end{exact}
\begin{flushright}
$\diamond$
\end{flushright}
\noindent In fact, a bijection $\nu$ between permutations in
$S_n(3\textrm{-}1\textrm{-}2,321)$ and Motzkin paths of length $n$
can be obtained as the composition of the map $K$ with the
well-known bijection between Dyck $n$-paths with no $DDD$ and
Motzkin $n$-paths, defined by replacing each $UDD$ with $D$ and
each remaining $UD$ with a horizontal step $H$ (see e.g.
\cite{clae}).

\begin{figure}[ht]
\begin{center}
\includegraphics[bb=98 635 263 698,width=.55\textwidth]{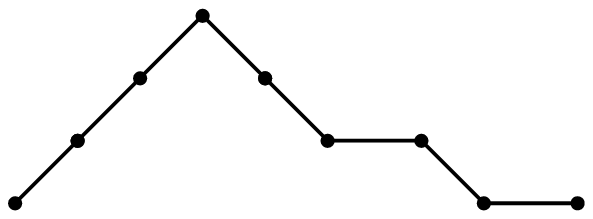}
\caption{The Motzkin path $\nu(4\,3\,5\,2\,6\,7\,1\,8)$}
\end{center}
\end{figure}

\subsection{The pattern $231$}

\noindent First of all, note that occurrences of the consecutive
pattern $231$ in the $3\textrm{-}1\textrm{-}2$ avoiding
permutation
$$\sigma=m_1\,w_1\,m_2\,w_2\,\ldots\,m_k\,w_k$$
correspond bijectively to those indices $i$ such that:
\begin{itemize}
\item[1.] $m_{i+1}=m_i+1$;
\item[2.] $w_{i+1}$ is nonempty.
\end{itemize}
In fact, consider the subword $m_i\,w_i\,m_{i+1}\,w_{i+1}$ and its
subword $a\, m_{i+1}\, b$, where $a$ is the rightmost symbol in
$w_i$ (or $m_i$ if $w_i$ is empty) and $b$ is the leftmost symbol
in $w_{i+1}$. Then, $a\, m_{i+1}\, b$ is order isomorphic to $231$
if and only if $a>b$. This happens whenever the two conditions
above hold. Each one of these occurrences corresponds to an
occurrence of $DUDD$ in the Dyck path $K(\sigma)$. The
distribution of $DUDD$ on Dyck paths was deeply studied (see
\cite{mansou}, \cite{sapo}, and \cite{oeis} seq.$\,$A116424). In
\cite{sapo}, the author deduces a functional equation satisfied by
the bivariate generating function of this distribution. Then:
\newtheorem{deduci}[yama]{Theorem}
\begin{deduci} We have
$$A^{231}(t,z)=\frac{1-(1-z)t^2-\sqrt{((1-z)t^2+1)^2-4t}}{2t(1-(1-z)t)}$$
that yields
$$a^{231}_{n,k}=\sum_{j=k}^{\left\lfloor\frac{n-1}{2}\right\rfloor}\frac{(-1)^{j-k}}{n-j}
{j\choose k}{n-j\choose j}{2n-3j\choose n-j+1}.$$
\end{deduci}
\begin{flushright}
$\diamond$
\end{flushright}
\noindent In particular:
\newtheorem{partic}[yama]{Proposition}
\begin{partic}
The number of permutations in $S_n$ that avoid both the pattern
$3\textrm{-}1\textrm{-}2$ and the consecutive pattern $231$ is
$$|S_n(3\textrm{-}1\textrm{-}2,231)|=\sum_{j=0}^{\left\lfloor\frac{n-1}{2}\right\rfloor}\frac{(-1)^{j}}{n-j}
{n-j\choose j}{2n-3j\choose n-j+1}.$$
\end{partic}

\subsection{The pattern $123$}

\noindent A $3\textrm{-}1\textrm{-}2$-avoiding permutation
$$\sigma=m_1\,w_1\,m_2\,w_2\,\ldots\,m_k\,w_k$$
contains an occurrence of the consecutive pattern $123$ if and
only if it contains two adjacent left-to-right maxima
$m_i\,m_{i+1}$, with $i>1$. In fact, in this case, $\sigma$
contains the $123$-subword $a\,m_i\,m_{i+1}$, where $a$ is the
symbol preceding $m_i$.

\noindent This consideration implies that occurrences of $123$ in
$\sigma$ correspond bijectively to occurrences of $DU^tDU$, with
$t>0$, in the Dyck path $K(\sigma)$.

\newtheorem{rogo}[yama]{Proposition}
\begin{rogo}\label{cioe}
The two statistics ''number of occurrences of $DDD$'' and ''number
of occurrences of $DU^tDU$'', $t>0$, are equidistributed on Dyck
$n$-paths.
\end{rogo}

\noindent \emph{Proof} Consider Deutsch's involution $\Delta$
described in Section \ref{bestia}. Denote by $f(P)$ the number of
occurrences of the subword $DDD$ in a given Dyck path $P$ and by
$g(P)$ the number of occurrences of $DU^tDU$, $t>0$, in $P$.

\begin{figure}[h]
\begin{center}
\includegraphics[bb=81 498 368 772,width=.6\textwidth]{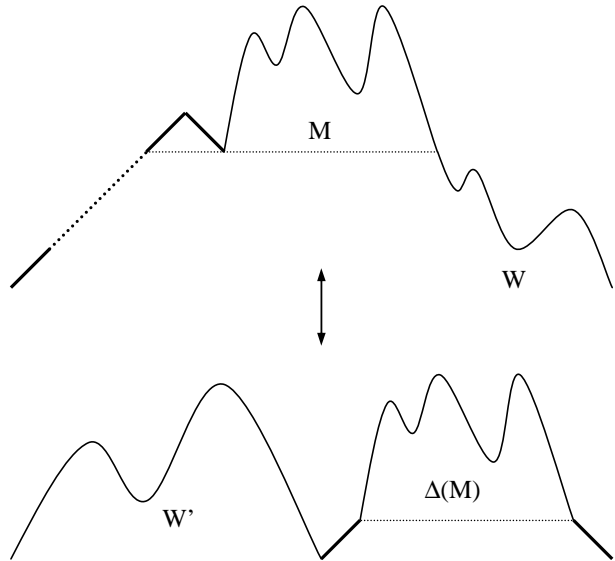}
\caption{The involution $\Delta$ maps a Dyck path starting with
$U^t\,D\,U$ to a Dyck path ending with two down
steps.}\label{prefa}
\end{center}
\end{figure}

\noindent  We prove that $g(P)=f(\Delta(P))$ by induction on the
semilength $n$ of the Dyck path $P$. The assertion is trivially
true for $n=0$. Fix an integer $n>0$ and assume that the assertion
holds for all Dyck paths of semilength less than $n$. Let $P$ be a
Dyck $n$-path, and $P=U\,A\,D\,B$ its first return decomposition.
It is easy to verify that:
$$f(P)=\left\{\begin{array}{ll}f(A)+f(B)+1&\textrm{if }A\textrm{ ends with }DD\\f(A)+f(B)&\textrm{otherwise}
\end{array}\right.$$
$$g(P)=\left\{\begin{array}{ll}g(A)+g(B)+1&\textrm{if }B\textrm{ begins with }U^tDU\\g(A)+g(B)&\textrm{otherwise}
\end{array}\right..$$
Since the semilengths of the two Dyck paths $A$ and $B$ are
strictly less than $n$, by induction hypothesis we have
$g(A)=f(\Delta(A))$ and $g(B)=f(\Delta(B))$. It is sufficient to
show that the involution $\Delta$ acts as follows: if a Dyck path
$P$ begins with a subword of type $U^t\,D\,U$, then $\Delta(P)$
ends with two down steps. Consider a Dyck path $P$ starting with
 $U^t\,D\,U$. Such a Dyck path can be decomposed into
$U^{t-1}\,U\,D\,M\,W$, where $M$ is an irreducible (and nonempty)
Dyck path. In this case, the Dyck path $\Delta(P)$ decomposes into
$\Delta(P)=W'\,U\,\Delta(M)\,D$ (see Figure \ref{prefa}). The path
$\Delta(M)$ is nonempty, hence $\Delta(P)$ must end with two down
steps, as desired.

\noindent Since $\Delta$ is an involution, we have also
$f(P)=g(\Delta(P))$. This proves that $\Delta$ maps every
occurrence of the subword $DU^tDU$ into an occurrence of the
subword $DDD$, and conversely.

\begin{flushright}
$\diamond$
\end{flushright}

\noindent As an immediate consequence, we have:

\newtheorem{dedu}[yama]{Theorem}
\begin{dedu}
The two generating functions $A^{123}(t,z)$ and $A^{321}(t,z)$
coincide.
\end{dedu}
\begin{flushright}
$\diamond$
\end{flushright}

\noindent Hence, by Proposition \ref{arrsub}, we can state the
following

\newtheorem{vedi}[yama]{Proposition}
\begin{vedi}
The number of permutations in $S_n$ that avoid both the pattern
$3\textrm{-}1\textrm{-}2$ and the consecutive pattern $123$ is
$$|S_n(3\textrm{-}1\textrm{-}2,123)|=M_n.$$
\end{vedi}
\begin{flushright}
$\diamond$
\end{flushright}
We submit that an explicit bijection between the set
$S_n(3\textrm{-}1\textrm{-}2,123)$ and the set of Motzkin
$n$-paths can be described as follows: let $\sigma$ be a
permutation in $S_n(3\textrm{-}1\textrm{-}2,123)$. The permutation
$\sigma$ is either of kind
$$\sigma=\sigma_1\ 1$$
or of kind
$$\sigma=\sigma_1\ 1\ t\ \sigma_2,$$
where $\sigma_1$ and $\sigma_2$ avoid $3\textrm{-}1\textrm{-}2$
and $123$ and $t=s+1$, where $s$ is the leftmost element appearing
in $\sigma_2$ (note that, in both cases, $\sigma_1$ and $\sigma_2$
can be the empty permutation). We define the bijection
$\mu:S(3\textrm{-}1\textrm{-}2,123)\to\mathscr{M}_n$ recursively
as follows:
\begin{itemize}
\item the empty permutation is mapped to the empty path;
\item if $\sigma=\sigma_1\ 1\ t\ \sigma_2$, then
$\mu(\sigma)=\mu(\sigma_1)\,U\,\mu(\sigma_2)\,D$;
\item if $\sigma=\sigma_1\ 1$, then
$\mu(\sigma)=\mu(\sigma_1)\,H$.
\end{itemize}

\noindent For example, consider the permutation
$\sigma=2\,4\,3\,1\,6\,5\,8\,7\in
S_8(3\textrm{-}1\textrm{-}2,123)$. Then,
$$\mu(\sigma)=\mu(2\,4\,3)U\mu(5\,8\,7)D=U\mu(3)DUU\mu(7)DD=UHDUUHDD.$$

\begin{figure}[h]
\begin{center}
\includegraphics[bb=97 631 268 676,width=.55\textwidth]{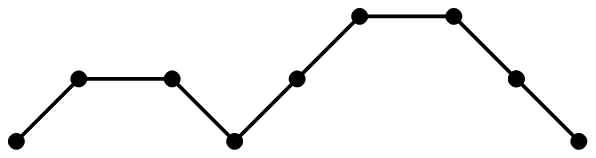}
\caption{The Motzkin path $\mu(2\,4\,3\,1\,6\,5\,8\,7)$.}
\end{center}
\end{figure}

\noindent We point out that the argumentations used in the proof
of Proposition \ref{cioe} can be iterated to get the following
more general result:

\newtheorem{mg}[yama]{Proposition}
\begin{mg}
The two bivariate generating functions $A^{12\cdots k}(t,z)$ and
$A^{k\cdots 21}(t,z)$ coincide.
\end{mg}

\subsection{The pattern $132$}

\noindent  A $3\textrm{-}1\textrm{-}2$-avoiding permutation
$$\sigma=m_1\,w_1\,m_2\,w_2\,\ldots\,m_k\,w_k$$
contains an occurrence of the consecutive pattern $132$ if and
only if there exists an index $i$, $i<k$, such that the subword
$m_i\,w_i\,m_{i+1}\,w_{i+1}$ verify the following conditions:
\begin{enumerate}
\item $m_{i+1}-m_i>1$;
\item $w_{i+1}$ is nonempty.
\end{enumerate}
In fact, in this case, $\sigma$ contains the $132$-subword
$a\,m_{i+1}\,b$, where $a$ is the symbol preceding $m_{i+1}$ in
$\sigma$ and $b=m_{i+1}-1$.

\noindent These remarks imply that occurrences of $132$ in
$\sigma$ correspond bijectively to occurrences of $DU^tDD$, with
$t>1$, in the Dyck path $K(\sigma)$.

\newtheorem{piange}[yama]{Proposition}
\begin{piange}\label{ilresto}
The two statistics ''number of occurrences of $DUDD$'' and
''number of occurrences of $DU^tDD$'', $t>1$, are equidistributed
on Dyck $n$-paths.
\end{piange}

\noindent \emph{Proof} Denote by $h(P)$ the number of occurrences
of the subword $DUDD$ in a given Dyck path $P$ and by $l(P)$ the
number of occurrences of $DU^tDD$, $t>1$, in $P$.

\begin{figure}[ht]
\begin{center}
\includegraphics[bb=90 506 342 734,width=.6\textwidth]{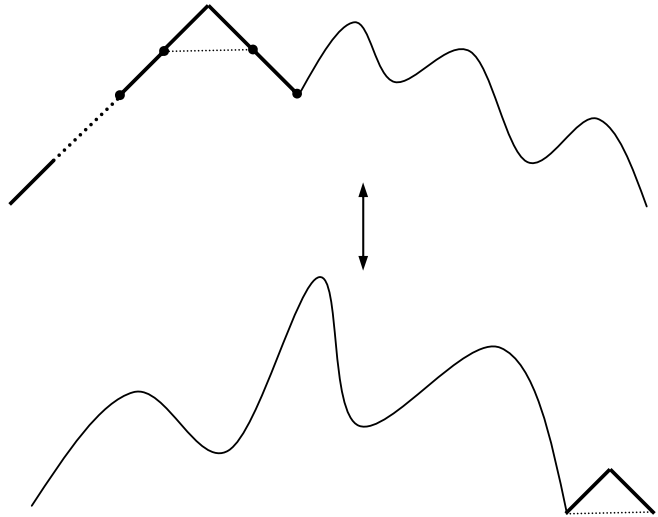}
\caption{The involution $\Delta$ maps a Dyck path starting with
$U^t\,D\,D$, $t>1$, to a Dyck path ending with
$U\,D$.}\label{fapre}
\end{center}
\end{figure}

\noindent We prove that $l(P)=h(\Delta(P))$, where $\Delta$ is
Deutsch's involution, by induction on the semilength $n$ of the
Dyck path $P$. The assertion is trivially true for $n=0$. Fix an
integer $n>0$ and assume that the assertion holds for all Dyck
paths of semilength less than $n$. Let $P$ be a Dyck $n$-path, and
$P=U\,A\,D\,B$ its first return decomposition. It is easy to
verify that:
$$h(P)=\left\{\begin{array}{ll}h(A)+h(B)+1&\textrm{if }A\textrm{ ends with }UD\\h(A)+h(B)&\textrm{otherwise}
\end{array}\right.$$
$$l(P)=\left\{\begin{array}{ll}l(A)+l(B)+1&\textrm{if }B\textrm{ begins with }U^tDD\ ,t>1\\l(A)+l(B)&\textrm{otherwise}
\end{array}\right.$$
Since the semilengths of the two Dyck paths $A$ and $B$ are
strictly less than $n$, by induction hypothesis we have
$l(A)=h(\Delta(A))$ and $l(B)=h(\Delta(B))$. It is sufficient to
show that the involution $\Delta$ maps a Dyck path $P$ that begins
with a subword of type $U^t\,D\,D$, $t>1$, to a path $\Delta(P)$
ending with $UD$. Consider a Dyck path $P$ starting with the
subword $U^t\,D\,D$, $t>1$. In this case, the last step of the
recursive procedure defining the map $\Delta$ maps the first peak
$UD$ of $P$ to the last irreducible component of the Dyck path
$\Delta(P)$ (see Figure \ref{fapre}). Since $\Delta(UD)=UD$, the
Dyck path $\Delta(P)$ ends with $UD$, as desired.

\noindent Recalling that $\Delta$ is an involution, we have also
$h(P)=l(\Delta(P))$. Hence, $\Delta$ maps every occurrence of the
subword $DU^tDD$, $t>1$, into an occurrence of the subword $DUDD$,
and viceversa.
\begin{flushright}
$\diamond$
\end{flushright}

\noindent Hence, we have:

\newtheorem{piuya}[yama]{Theorem}
\begin{piuya}
The two generating functions $A^{132}(t,z)$ and $A^{231}(t,z)$
coincide.
\end{piuya}

\begin{flushright}
$\diamond$
\end{flushright}

\section{Joint distributions}

\noindent The bijection $\Delta$ on Dyck paths induces an
involution $\hat{\Delta}=K^{-1}\circ\Delta\circ K$ on the set
$S_n(3\textrm{-}1\textrm{-}2)$. In particular, $\hat{\Delta}$ acts
on $S_3(3\textrm{-}1\textrm{-}2)$ as follows:
$$321 \longleftrightarrow^{\hspace{-.5cm}\Delta}\hspace{.2 cm} 123$$
$$231 \longleftrightarrow^{\hspace{-.5cm}\Delta}\hspace{.2 cm} 132$$
$$213 \longleftrightarrow^{\hspace{-.5cm}\Delta}\hspace{.2 cm} 213$$
In this section we prove that the action of $\hat{\Delta}$ on
$S_3(3\textrm{-}1\textrm{-}2)$ reveals to be paradigmatic for the
general case: the involution $\hat{\Delta}$ maps an occurrence of
a consecutive pattern $\tau\in S_3$ to an occurrence of the
consecutive pattern $\hat{\Delta}(\tau)$.\\

\noindent First of all, the statements of Propositions \ref{cioe}
and \ref{ilresto} can be reformulated in terms of $\hat{\Delta}$
as follows:

\newtheorem{delle}[yama]{Proposition}
\begin{delle}
Let $\sigma\in S_n(3\textrm{-}1\textrm{-}2)$. Then:
\begin{itemize}
\item[1.] $\sigma$ contains $k$ occurrences of the consecutive
pattern $123$ $\iff$ $\hat{\Delta}(\sigma)$ contains $k$
occurrences of the consecutive pattern $321$;
\item[2.] $\sigma$ contains $k$ occurrences of the consecutive
pattern $132$ $\iff$ $\hat{\Delta}(\sigma)$ contains $k$
occurrences of the consecutive pattern $231$.
\end{itemize}
\end{delle}
\begin{flushright}
$\diamond$
\end{flushright}

\noindent We now determine the behavior of the map $\hat{\Delta}$
with respect to the distribution of the pattern $213$:

\newtheorem{nessun}[yama]{Proposition}
\begin{nessun}
The two permutations $\sigma$ and $\hat{\Delta}(\sigma)$ have the
same number of occurrences of the consecutive pattern $213$.
\end{nessun}

\noindent \emph{Proof} Let $\sigma\in
S_n(3\textrm{-}1\textrm{-}2)$. As remarked in Subsection
\ref{semi}, each occurrence of $213$ in $\sigma$ corresponds to an
occurrence of $DDU$ in the Dyck path $K(\sigma)$. Hence, it is
sufficient to prove that, for any $P$, the two Dyck paths $P$ and
$\Delta(P)$ have the same number of occurrences of $DDU$.

\noindent Denote by $r(P)$ the number of occurrences of the
subword $DDU$ in $P$. We prove that $r(P)=r(\Delta(P))$ by
induction on the semilength $n$ of the Dyck path $P$. The
assertion is trivially true for $n=0$. Fix an integer $n>0$ and
assume that the assertion holds for all Dyck paths of semilength
less than $n$. Let $P$ be a Dyck $n$-path, and $P=U\,A\,D\,B$ its
first return decomposition. It is easy to verify that:
$$r(P)=\left\{\begin{array}{ll}r(A)+r(B)+1&\textrm{if both }A\textrm{ and }B\textrm{ are nonempty}\\h(A)+h(B)&\textrm{otherwise}
\end{array}\right.$$
Since the semilengths of the two Dyck paths $A$ and $B$ are
strictly less than $n$, by induction hypothesis we have
$r(A)=r(\Delta(A))$ and $r(B)=r(\Delta(B))$. Noting that both $A$
and $B$ are nonempty if and only if $\Delta(A)$ and $\Delta(B)$
are nonempty, we get the assertion.

\begin{flushright}
$\diamond$
\end{flushright}

\noindent Let $occ_{\tau}(\sigma)$ be the number of occurrences of
the pattern $\tau$ in the permutation $\sigma$. The preceding
results can be restated as follows:

\newtheorem{joi}[yama]{Theorem}
\begin{joi}
The triplets of statistics
\begin{itemize}
\item $(occ_{321},occ_{132},occ_{213})$ and $(occ_{123},occ_{231},occ_{213})$
\item $(occ_{321},occ_{231},occ_{213})$ and $(occ_{123},occ_{132},occ_{213})$
\end{itemize}
are equidistributed on $S_n(3\textrm{-}1\textrm{-}2)$, namely,
$$\sum_{\sigma\in S_n(3\textrm{-}1\textrm{-}2)}x^{occ_{321}(\sigma)}y^{occ_{132}(\sigma)}z^{occ_{213}(\sigma)}=
\sum_{\sigma\in
S_n(3\textrm{-}1\textrm{-}2)}x^{occ_{123}(\sigma)}y^{occ_{231}(\sigma)}z^{occ_{213}(\sigma)},$$
$$\sum_{\sigma\in S_n(3\textrm{-}1\textrm{-}2)}x^{occ_{321}(\sigma)}y^{occ_{231}(\sigma)}z^{occ_{213}(\sigma)}=
\sum_{\sigma\in
S_n(3\textrm{-}1\textrm{-}2)}x^{occ_{123}(\sigma)}y^{occ_{132}(\sigma)}z^{occ_{213}(\sigma)}.$$
\end{joi}
\begin{flushright}
$\diamond$
\end{flushright}

\end{document}